\documentclass[smallextended]{svjour3}
\usepackage[utf8]{inputenc}
\usepackage{babel}
\usepackage{hyperref}
\usepackage{xcolor}
\usepackage{color}
\usepackage{graphicx}
\usepackage{subcaption}
\usepackage{amssymb,amsmath}
\usepackage{booktabs}
\usepackage[linesnumbered,lined,commentsnumbered]{algorithm2e}
\usepackage{pgfplots}
\usepackage{mathtools}
\usepackage{placeins}

\usepackage{setspace}

\DeclareMathOperator{\im}{Im}     

\newcommand{\al}{\alpha}           
\newcommand{\sg}{\sigma}           
\newcommand{\sse}{\subseteq}       

\newcommand{\bR}{\mathbb{R}}        
\newcommand{\bS}{\mathbb{S}}        
\newcommand{\bZ}{\mathbb{Z}}        

\newcommand{\sA}{\mathcal{A}}       
\newcommand{\sE}{\mathcal{E}}       
\newcommand{\sG}{\mathcal{G}}       
\newcommand{\sV}{\mathcal{V}}       

\renewcommand{\geq}{\geqslant}      
\renewcommand{\leq}{\leqslant}      
\renewcommand{\o}{\circ}            
\renewcommand{\:}{\colon}           

\newcommand{\lb}{\left\{}       
\newcommand{\rb}{\right\}}        

\newcommand{\word}[1]{\quad\text{#1}\quad} 
\newcommand{\xto}[1]{\xrightarrow{#1}}      

\begin{document}
\title{Estimation of first-order sensitivity indices based on symmetric reflected Vietoris-Rips complexes areas}

\author{Alberto J. Hern\'andez \and Maikol Sol\'is \and Ronald A.
	Z\'u\~niga-Rojas}

\institute{Alberto J. Hern\'andez \at Centro de Investigaci\'on en Matem\'atica
	Pura y Aplicada (CIMPA), Escuela de Matem\'atica, Universidad de Costa Rica,
	San Jos\'e, Costa Rica \email{albertojose.hernandez@ucr.ac.cr}\and %
	Maikol Sol\'is (Corresponding author) \at Centro de Investigaci\'on en Matem\'atica Pura y Aplicada
	(CIMPA), Escuela de Matem\'atica, Universidad de Costa Rica, San Jos\'e, Costa
	Rica \email{maikol.solis@ucr.ac.cr}\and %
	Ronald A. Z\'u\~niga-Rojas \at Centro de Investigaci\'on Matemáticas y
	Meta-Matem\'aticas (CIMM), Escuela de Matem\'atica, Universidad de Costa Rica,
	San Jos\'e, Costa Rica \email{ronald.zunigarojas@ucr.ac.cr} }
\maketitle


\begin{abstract}
  In this paper we estimate the first-order sensitivity index of random
  variables within a model by reconstructing the embedding manifold of a
  two-dimensional cloud point. The model assumed has $p$ predictors and a
  continuous outcome $Y$. Our method gauges the manifold through a Vietoris-Rips
  complex with a fixed radius for each variable. With this object, and using the
  area and its symmetric reflection, we can estimate an index of relevance for
  each predictor. The index reveals the geometric nature of the data points.
  Also, given the method used, we can decide whether a pair of non-correlated
  random variables have some structural pattern in their interaction.
\end{abstract}{}








\keywords{  Sensitivity analysis; topological data analysis; invariant; symmetric
  difference; symmetric reflection}


\subclass{49Q12; 55U10; 14J33}

\section{Introduction}\label{sec:introduction}

Topological data analysis (TDA) is a recent area of research  that studies the
topological invariants of data point clouds and applies them to specific
problems in statistics, machine learning or data visualization. TDA provides
new insights to statistical problems and helps discover hidden patterns
that classical methods could not detect. Novel examples for TDA are in
\cite{LiTopological2019} \cite{RybakkenDecoding2019},
\cite{Paluzo-HidalgoPhilological2020}. A recent survey on the topic could be
found in \cite{OtterRoadmap2017}.

The basic workflow described by \cite{GhristBarcodes2008} for TDA has three
steps: 1. Convert the set of data points into a family of simplicial complexes
(a topological manifold); 2. Identify the persistent homology (Betti numbers) of
the whole complex; 3. Encode the persistence homology into a barcode
(\cite{GhristBarcodes2008}), a persistence diagram
(\cite{Cohen-SteinerStability2007}), landscapes (\cite{BubenikStatistical2015}),
or other similar structures.

In this workflow, the first step simplifies the random nature of the data points
into a single structure which can be tackled as a whole. The second step uses
algebraic topology to separate all the relevant features of the data. In
particular we can identify the Betti numbers which are the quantity of connected
components, holes, voids, etc. Betti numbers represent an example of a
topological invariant of the manifold that TDA uses to identify the different
features of the data. Euler characteristic and cohomology groups are also
examples of topological invariants of the complex that can be used. Finally, we
can encode the Betti numbers into simpler structures like barcode, persistence
diagram or landscapes.

One problem in statistics, which is of our interest, is sensitivity analysis
(SA). It explores how much one variable impacts an output through a model. The
model could be known, unknown, a computer code or even a meta-model
(\cite{GratietMetamodelBased2015}). Classic methods are widely known, using
techniques like screening methods (\cite{CampolongoScreening2011}), Montecarlo
simulation (\cite{IshigamiImportance1990}), Sampling (
\cite{SobolEstimating2007}, \cite{SaltelliVariance2010}), Non-parametric curves
(\cite{SolisNonparametric2019}), among others.

In this paper, we will focus on the use of TDA to determine the sensitivity of a
variable into a model. Using the geometric shape of the data, we will identify
the most relevant variables of the model. The relevance could come as having
monotonic or semi-monotonic patterns, holes inside the data or anomalous
patterns. The non-important variables should have flat or noisy shapes. Given
our interest  on detecting the shape of the data, we can also distinguish
the structured noise from the unstructured one.

In our framework all the information is provided by the data. The underlying
function or process that generated the data is unknown. Our aim is to determine
which variable has the greatest impact with respect to the output. To achieve
this, we will analyze the data for each variable separately, transform it into a
simplicial complex and then build its Vietoris-Rips complex. In this way we will
construct an estimate of the shape that contains the geometric features for each
variable and its output.

We have to clarify that constructing the Vietoris-Rips for a fixed neighborhood
radius \(\varepsilon\) will not show all the features of the data (see p.64 on
\cite{GhristBarcodes2008}). The persistent homology arisees when taking into
consideration a range of possible values for the neighborhood radii which can be
detected in a barcode in the form of continuous lines representing persisting
features of the topology (\cite{CollinsBarcode2004}). However, it is possible to
take the shorter-medium lines to extract more information
(\cite{StolzPersistent2017}) or use the disconnected components to cluster the
data (\cite{ChazalPersistenceBased2013}). In recent years new techniques used to
automatically separate topological features from noise straight from the barcode
have been studied. (\cite{AtienzaPersistent2019}). In this work, the radius is
chosen empirically, observing the barcodes and selecting the most prominent
features of the data.

In \cite{HernandezGeometrical2019} we constructed a geometric goodness-of-fit
index by using a similar method as the one described above . We were concerned
with the relationship of the input variable with respect to the output one. This
previous work lead to define a geometric correlation index between the
variables. However, it could not determine the true influence of a variable
inside a model. The most exceptional examples are when the variables have a
zero-sum pattern. In this case, their algorithm failed to acknowledge the
irrelevance of such variables.

To extend those ideas, we will use the reflection symmetry of an object to
compare it with itself. In this manner the zero-sum patterns and their own
reflections must have an almost identical shape which could allow us to identify
anomalous patterns. The comparison can also help distinguish pure noise,
structured noised and relevant variables. The creation of these symmetries are
possible thanks to some affine transformations on the vector space containing
the object. The work of \cite{CombesAutomatic2008} and \cite{NagarDetecting2019}
use them to identify symmetry patterns in 2D and 3D.

Finally, we estimate the area of the shape and that of its reflection. We use
the symmetric difference as an invariant between these two objects. If both are
close, in a geometric congruence sense, their symmetric difference should be
low, as a percentage of the total area encompassed by the superposition of
figures. This behavior shows that the point set has a noisy nature and its
correlation with the output variable is negligible. Otherwise, a high symmetric
difference area reveals that the original object and their symmetric difference
have strong irregularities concluding a relevant impact of this variable with
respect to the output.

We normalize our geometric indices, between $0$ and $1$, to determine if one
variable is noisy or relevant. Those indices were made using the symmetric
difference between a Vietoris-Rips complex and its reflected sibling. We use the
symmetric difference between the two objects to determine their level of
dissimilarity. This technique allows us to separate pure noise, structured noise
and relevant variables.

The paper is structured as follows: Section~\ref{sec:preliminaries} introduces
the framework to the sensitivity analysis. In this section, we review the
concept and construction of a Vietoris-Rips complex and the notions of
persistent homology to create balanced complexes. Some theory on symmetries is
also presented in this part. Section~\ref{sec:methodology} is the main section
of the paper. Here we present the methodology used to create the sensitivity
geometrical indices. We explain the Vietoris-Rips construction, the use of the
symmetric reflection and the index estimation. In Section~\ref{sec:results} we
present some simulations using our own package called {\tt topsa}. We show how
our method performs and order the relevant variables in a model. Finally, in
Section~\ref{sec:conclusions} we present our discussion about the method and
possible further research work.

\section{Preliminaries}
\label{sec:preliminaries}

\subsection{Sensitivity Analysis}
\label{sec:sensitivity-analysis}
Let $(X_1, X_2, \ldots, X_{p}) \in \mathbb {R} ^ {p}$ for $p \geq 1$ and
$Y \in \mathbb {R}$ two random variables. Define the non-linear regression model
as

\begin{equation}
	\label{eq:regression_nonlinear}
	Y = f (X_1, X_2, \ldots, X_{p}) + \varepsilon.
\end{equation}

Here $\varepsilon$ is a random noise independent of $(X_1, X_2, \ldots, X_{p})$.
The unknown function $m: \mathbb {R} ^ {p} \mapsto \mathbb {R}$ describes the
conditional expectation of $Y$ given the $p$-tuple $(X_1, X_2, \ldots, X_{p})$.
Suppose as well that $(X_{i1}, X_{i2}, \ldots, X_{ip}, Y_i)$ for
$i = 1, \ldots, n$ is a size $n$ sample for the random vector
$(X_{1}, X_{2},\ldots, X_{p}, Y) $ .

One popular method to determine the  relevant
variables of the model is the variance-based global sensitivity analysis, this method was
proposed by Sobol~\cite{SobolSensitivity1993} in 1993 and is based on the ANOVA
decomposition. He proved that if $f$ is a squared integrable function then it
could be decomposed in the unitary cube as:
\begin{equation}
	\label{eq:sobol-descomposition}
	Y = f = f_{0} + \sum_{i} f_{i} + \sum_{\substack{ij\\ i\neq j}} f_{ij} +
	\cdots+ f_{12\cdots p}
\end{equation}
where each term is also square integrable over the domain. Also, for
each $i$ we have $f_i = f_i(X_i)$, $f_{ij} = f_{ij}(X_i, X_j)$ and so
on. This decomposition has $2^{p}$ terms and the first one, $f_0$, is
constant. The remaining terms are non-constant functions.
Sobol\cite{SobolSensitivity1993} also proved that this representation is unique
if each term has zero mean and the functions are pairwise orthogonal.
Equation~\eqref{eq:sobol-descomposition} can be interpreted as the
decomposition of the output variable $Y$ into its effects due to the
interaction with none, one or multiple variables. Taking expectation
in equation~\eqref{eq:sobol-descomposition} and simplifying the
expression we get:
\begin{align*}
	f_0                 & = \mathbb{E}[Y]                                 \\
	f_i(X_{i})          & = \mathbb{E}[Y\vert X_i] - \mathbb{E}[Y]        \\
	f_{ij}(X_{i},X_{j}) & = \mathbb{E}[Y\vert X_i, X_j] - f_i - f_j - f_0 \\ 
\end{align*}
and so on for all combinations of variables.

Once with this orthogonal decomposition, we measure the variance of
each element. The global-variance method estimates the regression
curves (surfaces) for each dimension, removing the effects due to
variables in lower dimensions. Then, it gauges the variance of each
curve (surface) normalized by the total variance in the model. For the
first and second-order effects the formulas are: 

\begin{align}
	\label{eq:sobol-indices}
	S_{i}  & = \frac{\mathrm{Cov}(f_{i}(X_{i}),Y)}{\mathrm{Var}(Y)}
	= \frac{\mathrm{Var}(\mathbb{E}[Y\vert X_{i}])}{\mathrm{Var}(Y)}       \\
	S_{ij} & = \frac{\mathrm{Cov}(f_{ij}(X_{i},X_{j}),Y)}{\mathrm{Var}(Y)}
	= \frac{\mathrm{Var}(\mathbb{E}[Y\vert X_{i}, X_{j}])}{\mathrm{Var}(Y)}
	- S_{i} - S_{j}. \notag{}
\end{align}
The formulas for the higher order terms are obtained recursively.



\subsection{Simplicial Homology Background}
\label{sec:preliminaries-Vietoris-Rips-complex}

For the purpose of this paper a topological object is either a connected surface
or a connected directed graph. Given a geometric object define a $0$-simplex as a
point, frequently called a \emph{vertex}. Since we deal with finite sets of
data, taking coordinates on the Euclidean Plane $\sE = \bR^2$, we denote a
$0$-simplex as a point $p_j = (x_j, y_j)$ for $j = 1, \ \dots,\ n$.

If we join two distinct $0$-simplices, $p_0,\ p_1$, by an oriented line segment,
we get a $1$-simplex called an {\it edge}:
$\overline{p_{0}p_{1}} = (p_{1}-p_{0})$.

Consider now three non-colinear points $p_0,\ p_1,\ p_2$ as $0$-simplices.
Together they form three $1$-simplices: $\overline{p_0 p_1}$,
$\overline{p_0 p_2}$ and
$\overline{p_1 p_2} = \overline{p_0 p_2} - \overline{p_0 p_1}$. This last
equation shows that only two of them are linearly independent and span the
third one. The union of these three edges form a triangular shape, a $2$-simplex
called a {\em face}, denoted as $\triangle(p_0 p_1 p_2)$ that contains all the
points enclosed between the edges:
\[
  \triangle(p_0 p_1 p_2) = 
  \triangle^{2} = 
  \lb
  p\in \bR^3\: p = \sum_{j=0}^{2}\lambda_j p_j,\ \sum_{j=0}^{2}\lambda_j = 1, \lambda_j \geq 0 
  \rb
  \subseteq \bR^{3}.
\]
The notation $\triangle^{2}$ allows us geometrically to realize the $2$-simplex
as a subspace $\triangle^{2}\subseteq \bR^{3}$ of dimension $2$; roughly
speaking, $\triangle^{2}$ is the convex hull of {\em three} affinely-independent
points $p_0,p_1,p_2$.

As well as $2$-simplices, if we consider four non-coplanar points, we can
construct a $3$-simplex called \emph{tetrahedron}. A generalization of dimension
$n$ will be a convex set in $\bR^n$ containing $\lb p_0, p_1, \dots, p_n \rb$ a
subset of $n+1$ distinct points that do not lie in the same hyperplane of
dimension $n-1$ or, equivalently, that the vectors
$\lb \overline{p_0 p_j} = p_j - p_0\rb, 0 < j\leq n $ are linearly independent.
In such a case, we are denoting the points $\lb p_0, p_1, \dots, p_n \rb$ as
vertices, and the usual notation would be $[p_0, p_j]$ for edges,
$[p_0, p_1, p_2]$ for faces, $[p_0, \dots, p_4]$ for tetrahedra, and
$[p_0, \dots, p_n]$ for $n$\emph{-simplices}. The standard $n$-simplex is
usually denoted
\[
  \triangle^{n} = \lb p\in \bR^{n+1}\: p = \sum_{j=0}^{n}\lambda_j p_j,\ \sum_{j=0}^{n}\lambda_j = 1, \lambda_j \geq 0 \rb \subseteq \bR^{n+1},
\]
as we did above. A few words about dimension: as well as we did for the triangle
$\triangle^2$ as a subspace of $\bR^3$, this standard notation allows to realize
the $n$-simplex as a subspace $\triangle^n \subseteq \bR^d$ of dimension $n$
(here \(d\geq n\)).

A $\bigtriangleup$-complex $X$ (or simplicial complex $X$) is the topological
quotient of a collection of disjoint simplices identifying some of their faces
via a family of linear homeomorphisms $\lb \sigma_{\al} \rb_{\al \in \sA}$ that
preserve the order of the vertices. We use these to identify the $n$-simplices:
$e_{\al}^n$. Using simple words, we can think about $X$ as a collection of
$n$-simplices such that if $\triangle^k \subseteq \triangle^n \in X$ then
$\triangle^k \in X$, calling $\triangle^k$ as a {\em face} of $\triangle^n$, and
the geometrical-dimensional notation gives us the chance to interpret the
identification above as the union of the simplices in $X$ so that those
simplices only intersect along the shared faces, and hence, the simplicial
complex $X$ may be embedded in $\bR^{d}$, when the maximal simplices in $X$ are
of dimension $n \leq d$:
\[
X = 
\left(
\bigcup_{k = 0}^{n}\triangle^{k}
\right)
\Big/_{\sim} 
\subseteq \bR^{d}.
\]

Now, we may define the simplicial homology groups of a $\bigtriangleup$-complex
$X$ as follows. Lets consider the free abelian group $\bigtriangleup_n(X)$ with
open $n$-simplices $e_{\al}^n \sse X$ as basis elements. The elements of this
group, known as \emph{chains}, look like linear combinations of the form
\begin{equation}
	c = \sum_{\al}n_{\al}e_{\al}^n\label{chain001}
\end{equation}
with integer coefficients $n_{\al}\in \bZ$. We also could write chains as linear
combinations of characteristic maps
\begin{equation}
	c = \sum_{\al}n_{\al}\sg_{\al}\label{chain002}
\end{equation}
where every $\sg_{\al}\: \bigtriangleup^n\to X$ is the corresponding
characteristic map of each $e_{\al}^n$, with image the closure of $e_{\al}^n$.
So, $c \in \bigtriangleup_n(X)$ is a finite collection of $n$-simplices in $X$
with integer multiplicities $n_{\al}$. The boundary of the $n$-simplex
$[p_0,\dots, p_n]$ consists of the various $(n-1)$-simplices
\[
	[p_0,\dots,\hat{p}_j,\dots, p_n] = [p_0,\dots,p_{j-1},p_{j+1},\dots, p_n].
\]
For chains, the boundary of $c = [p_0,\dots,p_n]$ is an oriented $(n-1)$-chain of the form
\begin{equation}
	\partial c = \sum_{j=0}^n (-1)^j [p_0,\dots,\hat{p}_j,\dots,
	p_n]\label{chain003}
\end{equation}
which is a linear combination of faces. This allows us to define the
\emph{boundary homomorphisms} for a general $\bigtriangleup$-complex $X$,
$\partial_{n}\: \bigtriangleup_n(X)\to \bigtriangleup_{n-1}(X)$ as follows:
\begin{equation}
	\partial_{n}(\sg_{\al}) =
	\sum_{j=0}^n(-1)^j\sg_{\al}\big|{{}_{[p_0,\dots,\hat{p}_j,\dots, p_n]}}.
\end{equation}

Hence, we get a sequence of homomorphisms of abelian groups

\begin{equation*}
	\dots \to
	\bigtriangleup_n(X)\xto{\partial_n}
	\bigtriangleup_{n-1}(X)\xto{\partial_{n-1}}
	\bigtriangleup_{n-2}(X)\to
	\dots \to
	\bigtriangleup_1(X)\xto{\partial_1}
	\bigtriangleup_{0}(X) \xto{\partial_0}  0
\end{equation*}
where $\partial \o \partial = \partial_{n} \o \partial_{n-1} = 0$
for all $n$.
This is usually known as a \emph{chain complex}. Since $\partial_{n}
	\o
	\partial_{n-1} = 0$, the $\im(\partial_n) \sse
	\ker(\partial_{n-1})$, and so, we
define the $n^{th}$ \emph{simplicial homology group} of $X$ as the
quotient
\begin{equation}
	\label{eq:def-homology-group}
	H_n^{\bigtriangleup}(X) = \frac{\ker(\partial_n)}{\im(\partial_{n+1})}.
\end{equation}
The elements of the kernel are known as \emph{cycles} and the elements of the
image are known as \emph{boundaries}.

Easy computations of sequences give us the simplicial homology of some examples:
the circle $X = \bS^1:$
$$
	H_0^{\bigtriangleup}(\bS^1) \cong \bZ,\quad H_1^{\bigtriangleup}(\bS^1) \cong
	\bZ,\quad H_n^{\bigtriangleup}(\bS^1) \cong 0\word{for} n \geq 2,
$$
and the torus $X = T \cong \bS^1\times \bS^1:$
$$
	H_0^{\bigtriangleup}(T) \cong \bZ,\quad H_1^{\bigtriangleup}(T) \cong \bZ
	\oplus \bZ,\quad H_2^{\bigtriangleup}(T) \cong \bZ,\quad
	H_n^{\bigtriangleup}(T) \cong 0\word{for} n \geq 3.
$$
In a very natural way, one can extend this process to define
singular homology
groups $H_n(X)$. This process, nevertheless, is not trivial but
natural. If $X$ is a $\bigtriangleup$-complex with finitely many $n$-simplices, then
$H_n(X)$
(and of course $H_n^{\bigtriangleup}(X)$) is finitely generated.

The $n^{th}${\em -Betti number} of $X$ is the number $b_{n}$ of summands isomorphic to
the additive group $\bZ$. The reader may see \cite{HatcherAlgebraic2000} for
details. We will be interested on $0$-th, $1$-st and $2$-nd Betti numbers, since
they represent the generators of the set of vertices, edges and faces
respectively.

For purposes of this research work, $2$-simplices constitute the building blocks
of our Vietoris-Rips complex.




\subsection{Vietoris-Rips Complex and Persistent Homology}
\label{sec:persistent-homology}


For the general case, let us consider a set of data points $D\subseteq \bR^{n}$,
and let $\varepsilon > 0$.

\begin{definition}
  The {\em Vietoris-Rips complex} of $D$ at scale $\varepsilon$ is defined as
    \[
      \sV_{\varepsilon}(D) := \lb \triangle^{k} \subseteq D:\ || p_{i} - p_{j} || \leq \varepsilon, \forall p_{i},p_{j} \in \triangle^{k} \rb
    \]
    where $||\cdot ||$ represents the Euclidean norm on $\bR^{n}$.
\end{definition}
From the last definition, we can see that the simplices in
$\sV_{\varepsilon}(D)$ have vertices that are at a distance less or equal than
$\varepsilon > 0$.

Several authors, among them Zomorodian~\cite{ZomorodianFast2010} as pioneer,
consider a fixed value $\varepsilon_{0}$ and compute $\sV_{\varepsilon_{0}}(D)$
to compute then the complex at any other scale $\varepsilon < \varepsilon_{0}$,
using a {\em weight function} $w:\ \sV_{\varepsilon_0}\to \bR_{+}$ defined as:
\[
w(\triangle^{k}) : = 
\left\{
\begin{array}{c r}
  ||p_{i} - p_{j}||
  & \textmd{if}\ \triangle^{k} = \lb p_{i},p_{j}\rb\\
  & \\
  \max \lb w(\triangle^{\ell}):\ \triangle^{\ell}\subseteq \triangle^{k} \rb & \textmd{otherwise.}
\end{array}
\right.
\]

At the end of the day, the weight value $w(\triangle^{k})$ will be the minimum $\varepsilon$ such that the simplex $\triangle^{k}$ enters the Vietoris-Rips-complex $\sV_{\varepsilon_0}(D)$.


\begin{definition}[Neighborhood graph]\label{def:neigborhood-graph}
    Given a set of data points $D$, and considering its non-oriented graph $\sG = V\cup E$ as the union of vertex $V = D$ and edges $E$, we define {\em the neighborhood graph} as the pair $(\sG, w)$ where $w:\ E\to \bR_{+}$ is the weight function defined on the edges.
\end{definition}

\begin{definition}[Vietoris-Rips neighborhood graph]\label{def:Vietoris-Rips-neigborhood-graph}
	Given $V \subseteq \mathbb{R}^{n}$ a set of vertices (data points $D = V$) and
  a scale parameter $\varepsilon\in \mathbb{R}$, the {\em Vietoris-Rips
    neighborhood graph} (Vietoris-Rips-neighborhood graph) is defined as
  $G_{\varepsilon}(V) = V\cup E_{\varepsilon}(V)$ where
	\begin{equation*}
		\label{eq:Vietoris-Rips-neigborhood} E_{\varepsilon}(V)  = \{\{p_{i},p_{j}\} \mid p_{i},p_{j}
		\in V \wedge p_{i}\neq p_{j}   \wedge ||p_{i} - p_{j}|| \leq \varepsilon \}.
	\end{equation*}
\end{definition}

From now on, for our purposes we will consider a geometric object to be the Vietoris-Rips-complex $\mathcal{V}$ of a given set of points or vertices $D = V \subseteq \mathbb{R}^2$.

\begin{definition}[Vietoris-Rips complex]\label{def:Vietoris-Rips-expansion}
	Given a set of vertices $V\subseteq \bR^{2}$ and its neighborhood graph
  $\sG_{\varepsilon}(V) = V\cup E_{\varepsilon}(V)$ for some $\varepsilon > 0$,
  their Vietoris-Rips complex $\mathcal{V} = \mathcal{V}(\sG_{\varepsilon})$ is
  defined as the union of $V$, $E_{\varepsilon}(V)$ and the set of 2-simplices
  $T$ that appear on $G_{\varepsilon}(V)$.
\end{definition}

\begin{definition}[Vietoris-Rips expansion]
	\label{def:Vietoris-Rips-expansion}
	Given a neighborhood graph $\sG_{\varepsilon}(V)$, their {\em Rips complex} $\mathcal{R}_{\varepsilon}$ is defined as all the edges of a simplex $\triangle^{k}$ that
	are in $\sG_{\varepsilon}(V)$. In this case $\triangle^{k}$ belongs to $\mathcal{R}_{\varepsilon}$. For $\sG_{\varepsilon}(V) = V\cup E_{\varepsilon}(V)$, we have
	\begin{equation*}
	\mathcal{R}_{\varepsilon} = V \cup E_{\varepsilon}(V) \cup \left\{\triangle^{k} \Big| {\binom{\triangle^{k}}{2}}
	\subseteq E_{\varepsilon}(V)\right\}.
	\end{equation*}
	where $\triangle^{k}$ is a simplex of $\sG_{\varepsilon}(V)$.
\end{definition}

Let us consider a parametrized family of spaces: a sequence of Rips complexes
$\{\mathcal{R}_j\}_{j=1}^{n}$ associated to a specific point cloud data for an
increasing sequence of radii $(\varepsilon_j)_{j=1}^{n}$. Instead of considering
the individual homology of the complex $\mathcal{R}_j$, consider the homology of
the inclusion maps
\[
  \mathcal{R}_1 \xrightarrow{i} \mathcal{R}_2 \xrightarrow{i} \dots \mathcal{R}_{n-1} \xrightarrow{i} \mathcal{R}_n
\]
{\em i.e.}, consider he homology of the iterated induced inclusions
\[
 i_{*}\colon 
 H_{*}(\mathcal{R}_j)
 \to
 H_{*}(\mathcal{R}_k)
 \quad 
 \textmd{for}
 \quad 
 j < k.
\]
These induced homology maps tell us which topological features persist. The
persistence concept says how Rips complexes become a good approximation to
\u{C}ech complexes (For full details see \cite{DeSilvaCoverage2007},
\cite{ZomorodianComputing2005}).

\begin{lemma}[Lemma~2.1.\cite{GhristBarcodes2008}]
 For any radio $\varepsilon > 0$ there are inclusions
 \[
  \mathcal{R}_{\varepsilon}
  \hookrightarrow
  \mathcal{C}_{\varepsilon \sqrt{2}}
  \hookrightarrow
  \mathcal{R}_{\varepsilon \sqrt{2}}.
 \]
\end{lemma}

To work with {\em persistent homology} in general, start considering a {\bf
  persistence complex} $\{C_{*}^{i}\}_{i}$ which is a sequence of chain
complexes joint with chain inclusion maps
$f^{i}\colon C_{*}^{i}\to C_{*}^{i+1}$, according to the tools we are working
with, we can take a sequence of Rips or \u{C}ech complexes of increasing radii
$(\varepsilon_j)_j$.

\begin{definition}
  For $j < k$, we define the $(j,k)$-persistent homology of the persistent
  complex $C = \{C_{*}^{i}\}_{i}$ as the image of the induced homomorphism
  $f_{*}\colon H_{*}(C_{*}^{j})\to H_{*}(C_{*}^{k})$. We denote it by
  $H^{j\to k}_{*}(C)$.
\end{definition}
 
\begin{definition}[Barcodes]
  The {\em barcodes} are graphical representations of the persistent homology of
  a complex $H^{j\to k}_{*}(C)$ as horizontal lines in a $XY$-plane whose
  corresponds to the increasing radii $\{\varepsilon_j\}_j$, and whose vertical
  $Y$-axis corresponds to an ordering of homology generators. Roughly speaking,
  a barcode is the persistent topology analogue to a Betti number.
\end{definition}

\begin{figure}
	\centering
	\includegraphics[width=\linewidth]{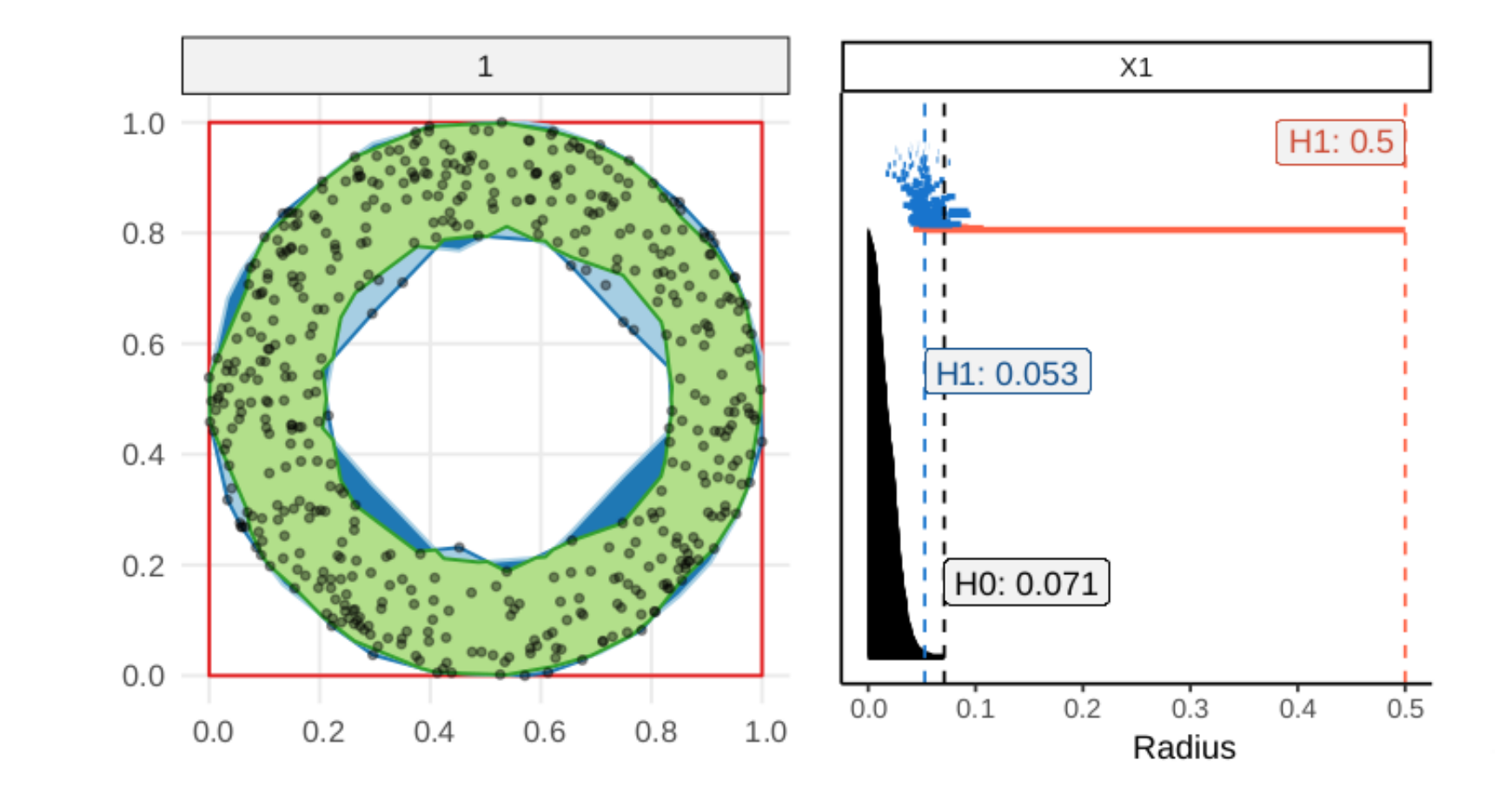}
	\caption{Manifold representing a circle with one hole and its corresponding barcode}
	\label{fig:circ1-ex-barcode}
\end{figure}

Figure~\ref{fig:circ1-ex-barcode} represents the barcodes for one set of data points arranged in a circle and having one hole in the middle (a doughnut). The central hole persists along multiple radius until it disappeared with \(\varepsilon = 0.5\). The long red line represent this feature persistence on $H_1$. The blue short lines are the noisy hole that appear and disappear briefly on the object. The black lines on the bottom are the number of connected components ($H_0$). 

\subsection{Symmetries and Affine Geometry}
\label{sec:symm-affine-geom}

Considering the Euclidean Plane $\mathbb{R}^{2}$, there are different
geometries we can get by considering the action of different groups on
it. Each group preserving a different set of properties that we would
like to study. For example, for a point $p\in \mathbb{R}^{2}$,
consider the action $\phi(p) = Ap + b$ where $A$ is a $2\times 2$
unitary matrix and $b \in \mathbb{R}^{2}$, this action preserves
distance, a concept central to the Euclidian Geometry in a classical
sense. Of course, the set of transformations described above is a
group under composition and congruence is an invariant under the
action of this set of transformations.

If we allow our set of transformations to be of the form $\phi(p) = Ap
	+ b$ with $A$ an invertible matrix, then the geometry we get is more
relaxed and the set of invariants of this action varies. This geometry
is called \textit{Affine Geometry}, and since not all invertible
matrices are unitary, affine geometry does not leave congruence in the
Euclidian sense as an invariant. Of course all isometries are affine
transformations.

Since affine transformations preserve both parallelism and ratios
along a given line, then it preserves ratios along parallel lines
which implies that it preserves axial symmetry and that it carries
within it other related properties in point sets.

In this study we will take advantage of this properties of affine
geometry to talk about certain geometric qualities of point sets that
are fundamental to understand sensitivity.

\section{Methodology}\label{sec:methodology}

In \cite{HernandezGeometrical2019}   it is defined the way that the Vietoris Rips
is constructed to determine the homological structure of the data. The
method is based in the work of \cite{ZomorodianFast2010} on the use of
cliques to dertermine the 2-simplex of the data.

In this work we determined a geometric correlation between the input
and output variables. The algorithm allows us to measure the geometric
pattern present in the data. The method consisted in three steps: 1.
Build the Vietoris-Rips complex for a single given radius; 2. Estimate the
bounding box contained inside the Vietoris-Rips complex; and 3. Compare the areas
of the bounding box and the Vietoris-Rips complex to determine the percentage of
blank space remained in the box.

The mentioned method is incapable to notice if one variable is
relevant or not to the model. In other words, a variable with a
non-noisy structure could be meaningless to explain the output. The
behavior is explained as it does not matter what value takes the input
variables, it does not get reflected over the output one. Even if
there is some structural pattern, those values nullify each other to
produce a non-influential process.

To explore further the intrinsic structure of the Vietoris-Rips complex, take the
classic Ishigami model $f(X_1, X_2, X_3) = \sin X_1 + 7 \sin^2 X_2 +
	0.1 X^3_3 \sin X_1$ where $X_d \sim \text{Uniform}[-\pi, \pi]$  and $d
	= 1, 2, 3$.  The reader can find a complete discussion of this model
in \cite{SaltelliScreening2009}.

\begin{figure}
	\centering \includegraphics[width=0.8\textwidth]{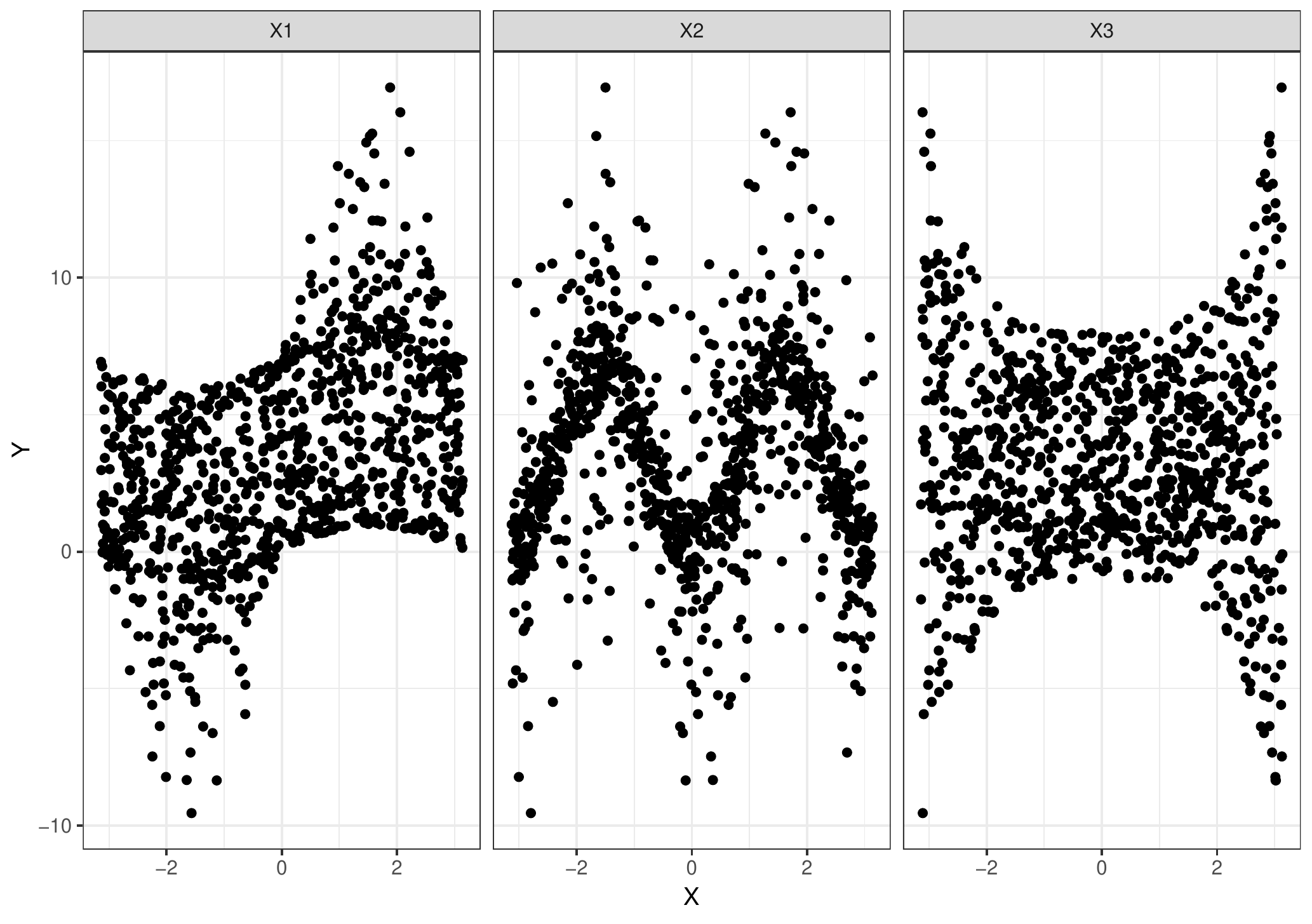}
	\caption{The Ishigami model plot with $n=1000$ random points.
		Theoretically, $S_{1} = 0.32$, $S_{2}=0.44$, $S_{3} = 0$ }
	\label{fig:ishigami-model}
\end{figure}

Notice from Figure~\ref{fig:ishigami-model} that the three variables have a
clear pattern. The variable $X_{1}$ is increasing, variable $X_{2}$ has a
\emph{M}-shape and variable $X_{3}$ has bow-tie shape. In the three cases,
\cite{HernandezGeometrical2019} estimated that the percentage of blank space is
around 50\%.

Nevertheless, the Sobol sensitivity indexes declare that the most important
variable is $X_{2}$, followed by $X_{1}$ and finally $X_{3}$ has a theoretical
null influence. One outstanding characteristic is how behaves the cloud of
points around the line middle line $\overline{Y} = (\min{Y_i} + \max{Y_i})/2$.
For the first two variables, we notice how the points upward the line does not
produce the same pattern on the downwards sections. However, for the third
variable $X_{3}$, both sections behave similar simulating a \emph{mirror
  effect}.

Thus, to measure the \emph{mirror effect} on a model, we propose the
following procedure:

\begin{enumerate}
	\item For each variable $X_{k}$, $k=1,\ldots,p$ and $Y$; project them
	      in the $\mathbb{R}^{2}$ plane.
	\item Estimate $(X_{\min}, Y_{\min} ) = (\min \{X_{ki}\}, \min
		      \{Y_{i}\})$ for  $i=1\ldots,n$.
	\item Recenter all the points $\tilde{X}_{ki} = X_{ki} - X_{\min}$ and
    $\tilde{Y}_{ki} = Y_{i} - Y_{\min}$.

	\item Create the Vietoris-Rips-complex $\mathcal{V}_k$ according to the procedure
	      in  \cite{HernandezGeometrical2019} for the cloud of points formed with
	      $(\tilde{X}_{ki},\tilde{Y}_{i})$.
	\item Create the symmetrical reflection $\mathcal{V}_k^{\text{SR}}$ of
	      the Vietoris-Rips-complex and set it in the same plane of $\mathcal{V}_k$ (see
	      Section~\ref{sec:affine-transf-vr}).
  \item Estimate the symmetric difference between both objects:
    $\mathcal{V}_k^{\triangle} = \mathcal{V}_k \triangle \mathcal{V}_k^{\text{SR}}$.
    Calculate the Area of $\mathcal{V}_k^{\triangle}$ and denote it as
    $\mathrm{Area}(\mathcal{V}_{\triangle})$. With those areas estimate the
    sensitivity index in Section~\ref{sec:estimation-sensitivity-index},
	      \begin{equation*}
		      S_{i}^{\text{Geom}} =
		      \frac{\mathrm{Area}(\mathcal{V}_k^{\triangle})}{2\
			      \mathrm{Area}(\mathcal{V}_k)}
	      \end{equation*}

\end{enumerate}

In the next sections we will describe in detail these steps.

\subsection{Affine transformation of the Vietoris-Rips
	complex}\label{sec:affine-transf-vr}

Recall that the Vietoris-Rips-complex is the union of the i-th complexes for $i
	\in \{0,1,2,\ldots\}$ of a datapoint cloud. The 2-simplexes are formed
with 3 vertex $(p_{k\alpha_{s}}, p_{k\beta_{s}}, p_{k\gamma_{s}})$
where the subindices $\alpha_{s}$, $\beta_{s}$ and $\gamma_{s}$ are
permutations with repetitions from the sequence $1,\ldots, n$ and
$s=1,\ldots, b_{2}$, where $b_{2}$ represents the second Betti number. Recall that each point has
associated a pairwise coordinates, e.g.,
$p_{k\alpha_{i}}=(X_{k\alpha_{s}},Y_{\alpha_{s}})^{\top}$.

Define the affine transformation $\varphi(p) = Ap +b$ where

\begin{align*}
	p & =(x,y)^{\top} \in \mathbb{R}^{2}         \\
	A & = \begin{pmatrix}
		1 & 0  \\
		0 & -1
	\end{pmatrix}             \\
	b & = (0, 2\ Y_{\operatorname{mid}} )^{\top}
\end{align*}

Applying the affine transformation $\varphi$ to each 2-simplex, we can
create the Vietoris-Rips-complex reflection over the original one. However, this
procedure fails if the data crosses the $x$-axis due to the reflection
get distorted.

To remedy this, we first have to shift the datapoints to obtain a
clear reflection through $\varphi$. To achieve this, recenter all the
points to $\tilde{X}_{ki} = X_{ki} - X_{\min}$ and $\tilde{Y}_{ki}
	=
	Y_{i} -
	Y_{\min}$, this way all the data will be contained in the first
quadrant.

To simplify the notation for the rest of the paper, we will assume
that every datapoint has been shifted according the last paragraph.

Thus, for a single 2-simplex $[p_{1},p_{2},p_{3}]$ defined in
Section~\ref{sec:preliminaries-Vietoris-Rips-complex} we can transform it
applying the affine transformation explained in
Section~\ref{sec:symm-affine-geom}. Therefore, we have

\begin{equation}
	\label{eq:phi-2-simplex}
	\varphi([p_{1},p_{2},p_{3}]) = \lambda_{1} \varphi(p_{1}) + \lambda_{2}
	\varphi(p_{2}) + \lambda_{3} \varphi(p_{3}).
\end{equation}

where $\lambda_{1} + \lambda_{2} + \lambda_{3} = 1$. %
Here we have abused of notation in $\varphi([p_{1},p_{2},p_{3}])$. It
means, applying the transformation $\varphi$ to the whole 2-simplex
$[p_1, p_2, p_3]$. By linearity, we apply $\varphi$ to inner point of
$[p_1,p_2,p_3]$, therefore it turns into the right side of
Equation~\eqref{eq:phi-2-simplex}.

Developing the last equation for a single point $p_{s} =
	(x_{s},y_{s})$ we have,

\begin{align*}
	\varphi(p_{s}) & =  Ap_{s} +b                 \\
	               & =
	\begin{pmatrix}
		1 & 0  \\
		0 & -1
	\end{pmatrix}
	\begin{pmatrix}
		x_{s} \\
		y_{s}
	\end{pmatrix} +
	\begin{pmatrix}
		0                          \\
		2   Y_{\operatorname{mid}}
	\end{pmatrix} \\
	               & = \begin{pmatrix}
		x_{s}                             \\
		2   Y_{\operatorname{mid}} -y_{s}
	\end{pmatrix} \\
\end{align*}

Therefore, gathering this results for the three points in a 2-simplex
$T$,

\begin{align*}
	\varphi([p_{1},p_{2},p_{3}]) & = \lambda_{1}
	\begin{pmatrix}
		x_{1}                               \\
		2\   Y_{\operatorname{mid}} - y_{1}
	\end{pmatrix}
	+ \lambda_{2}  \begin{pmatrix}
		x_{2}                               \\
		2\   Y_{\operatorname{mid}} - y_{2}
	\end{pmatrix}
	+  \lambda_{3y} \begin{pmatrix}
		x_{3}                               \\
		2\   Y_{\operatorname{mid}} - y_{3}
	\end{pmatrix}
\end{align*}

Finally, we define the symmetric reflection of $\mathcal{V}(G_{\varepsilon})$ as

\begin{equation*}
	\mathcal{V}_{k}^{\text{SR}}(G_{\varepsilon}) =
	\varphi(\mathcal{V}_{k}(G_{\varepsilon})) = \bigcup_{s=1}^{b_{2}}
	\left\{\varphi([p_{k\alpha_{s}},p_{k\beta_{s}},p_{k\gamma_{s}}])\right\}
\end{equation*}
where $b_{2}$ represents the second Betti number.

Here we apply the affine transformation $\varphi$ to all the 2-simplex
of the complex, and the join them into another complex.

For the rest of the paper, we will denote $\mathcal{V}_{k}\coloneqq
	\mathcal{V}_{k}(G_{\varepsilon})$ and $\mathcal{V}_{k}^{\text{SR}}\coloneqq
	\mathcal{V}_{k}^{\text{SR}}(G_{\varepsilon})$

\subsection{Estimation of the geometrical sensitivity index}
\label{sec:estimation-sensitivity-index}

Once the symmetric reflection $\mathcal{V}_{k}^{\text{SR}}$ is
estimated, we have to use it to determine if the variable $X_{k}$ is
relevant or  not. We mentioned before that one characteristic for the
non-relevant variables is their symmetry through the $y$-middle axis.
For this purpose, we will overlap the complexes
$\mathcal{V}_{k}^{\text{SR}}$ with respect to $\mathcal{V}_{k}$ and
gauge how much they are equal.

The symmetric difference between both complexes give the amount of area in
either one but not in the intersection. We will denote the new set as

\begin{equation*}
	\mathcal{V}_{k}^{\triangle} = \left(\mathcal{V}_{k} \cup
		\mathcal{V}_{k}^{\text{SR}}\right) \setminus   \left(\mathcal{V}_{k} \cap
		\mathcal{V}_{k}^{\text{SR}}\right).
\end{equation*}

Therefore, if the area of $\mathcal{V}_{k}^{\triangle}$ is small, it
means that the union of original complex and the transformed one, is
almost the same that their intersection. Otherwise, the points
produced a particular pattern that the symmetric reflection cannot
imitate. However, the complex $\mathcal{V}_{k}^{\text{SR}}$ replicate
the same behavior of  $\mathcal{V}_{k}$ in the opposite direction.

Once with the objects $\mathcal{V}_{k}$  and
$\mathcal{V}_{k}^{\triangle}$ lets define the geometric sensitivity
index as

\begin{equation}
	\label{eq:geometric-sensitivity-index} S_{i}^{\text{Geom}} =
	\frac{\mathrm{Area}(\mathcal{V}_k^{\triangle})}{2\
		\mathrm{Area}(\mathcal{V}_k)}
\end{equation}


The index $S_{i}^{\text{Geom}}$ measures how much the
$\mathcal{V}_{k}$ and $\mathcal{V}^{\text{SR}}_{k}$ are asymmetrical.

We can identify two cases:

\begin{enumerate}
	\item If $X_k$ is irrelevant to model $f$
	      (Equation~\eqref{eq:regression_nonlinear}) a sensitivity index should
	      be 0. The variable $X_k$ acts as a random input to the model causing
	      a random pattern for the output.

	      Here, we have that    $\mathcal{V}_{k}$ and
	      $\mathcal{V}^\star_{k}$ share the same structure (almost
	      everywhere) due to the random interaction between $X_k$ and
	      $Y$.

	\item If $X_k$ is important to the model, thus its values cause a
	      defined pattern to the output $Y$. Those patterns, in average, could
	      be increasing, decreasing or a combination of both, for instance.

	      In any case, the complexes $\mathcal{V}_{k}$ and
	      $\mathcal{V}^\star_{k}$ have a non-negligible difference.

\end{enumerate}

Considering these two cases, our index measures how asymmetrical are
$\mathcal{V}_{k}$ and $\mathcal{V}^\star_{k}$ through the
$\mathrm{Area}(\mathcal{V}_k^{\triangle})$. If the variable is
irrelevant or random $\mathcal{V}_{k} \cong \mathcal{V}^\star_{k}$,
which implies $\mathcal{V}_k^{\triangle} \cong \emptyset$. Therefore,
the area  of $\mathcal{V}_k^{\triangle}$ result in a negligible value.

In the other case,  $\mathcal{V}_{k} \ncong \mathcal{V}^\star_{k}$. In
the extreme case where the intersection of both sets are of measure
zero, we have
\begin{equation*}
	\mathrm{Area}(\mathcal{V}_k^{\triangle}) =
	\mathrm{Area}(\mathcal{V}_k) + \mathrm{Area}(\mathcal{V}^\star_k) = 2\
	\mathrm{Area}(\mathcal{V}_k).
\end{equation*}

Equation~\eqref{eq:geometric-sensitivity-index} summarizes both cases
to get a index between 0 to 1.

\section{Results}
\label{sec:results}

The methodology in this paper was compiled in a \textit{R}
(\cite{RCoreTeamLanguage2020}) package called  \emph{ \texttt{topsa}: Topological
	Sensitivity Analysis} (\url{https://cran.r-project.org/package=topsa}.

The additional packages used were  \texttt{TDA} (\cite{FasyTDA2019})
for estimating the Vietori-Rips structure, the package \texttt{sf}
(\cite{PebesmaSimple2018}) for all the geometric estimations and
the package \texttt{ggplot2} (\cite{WickhamGgplot22009}) for all the graphic
outputs.

For all the settings, we sample $n=1000$ points with the distribution
specified in each case. Due to the number of points, we choose the
quantile $5\%$ for each case to determine the radius of the
neighborhood graph. Further insights about this choosing will be
presented in the conclusions section

We will consider five settings, each one with different topological
features. The cases are not exhaustive and there are other settings
with interesting features as well. However, through those examples we
could show the capabilities of our models to detect sensitivity in
each variable.

\subsection{Theoretical examples}\label{ssec:theoretical-examples}
The examples considered are the following:

\paragraph{Linear} This is a simple setting with
\begin{equation*}
	Y = 2 X_{1} + X_{2}
\end{equation*}
and $X_{3}$ is an independent random variable. We set
$X_{i}\sim\mathrm{Uniform}(-1,1)$ for $i=1,\ldots, 3$.


\paragraph{Circle with hole:} The model in this case is
\begin{equation*}
	\begin{cases}
		X_{1} = r\cos(\theta) \\
		Y = r\sin(\theta)
	\end{cases}
\end{equation*}
with $ \theta \sim \mathrm{Uniform} (0,2\pi)$ and $r \sim \mathrm{Uniform}
	(0.5,1)$. This form creates a circle with a hole in the middle.

\paragraph{Connected circles with holes:} The model consists in three different circles, where we set $ \theta \sim \mathrm{Uniform} (0,2\pi)$:
\begin{enumerate}
	\item Circle centered at $(0,0)$ with radius between 1.5 and 2.5:
	      \begin{equation*}
		      \begin{cases}
			      X_{1} = r_1\cos(\theta) \\
			      Y = r_1\sin(\theta)
		      \end{cases}
	      \end{equation*}
	      where $r_{1}\sim \mathrm{Uniform} (1.5,2.5)$.
	\item Circle centered at $(3.5,3.5)$ with radius between 0.5 and 1:
	      \begin{equation*}
		      \begin{cases}
			      X_{2} - 3.5 = r_2\cos(\theta) \\
			      Y -3.5 = r_2\sin(\theta)
		      \end{cases}
	      \end{equation*}
	      where $r_{2}\sim \mathrm{Uniform} (0.5,1)$.
	\item Circle centered at $(-4,4)$ with radius between 1 and 2:
	      \begin{equation*}
		      \begin{cases}
			      X_{3} + 4 = r_3\cos(\theta) \\
			      Y - 4 = r_3\sin(\theta)
		      \end{cases}
	      \end{equation*}
	      where $r_{3}\sim \mathrm{Uniform} (1,2)$.
\end{enumerate}

\paragraph{Ishigami:} The final model is
\begin{equation*}
	Y = \sin X_1 + 7\ \sin^2 X_2 + 0.1\ X_3^4 \sin X_1
\end{equation*}
where $X_i\sim \mathrm{Uniform}(-\pi, \pi)$ for $i=1,2,3$, $a = 7$ and $b =
	0.1$.

\subsection{Numerical results}

The results in this section were estimated with the code \texttt{topsa}. The
figures are estimated Vietoris-Rips-complex for each input $X_{i}$ with respect
to the output variable $Y$. Also, we represent the symmetric
reflection of this complex and the symmetric difference between them.
The table below each figure presents the radius used to build the
neighborhood graph, the estimated areas of the manifold object
($\mathrm{Area}(\mathcal{V})$), the area of the reference box
contained the complex ($\mathrm{Area}(B)$), the geometrical
correlation index estimated in \cite{HernandezGeometrical2019}
($\rho^\mathrm{Geom}$), and the geometric sensitivity index
($S^\mathrm{Geom}$).

\begin{figure}[ht!]
	\centering
	\includegraphics[width=\textwidth]{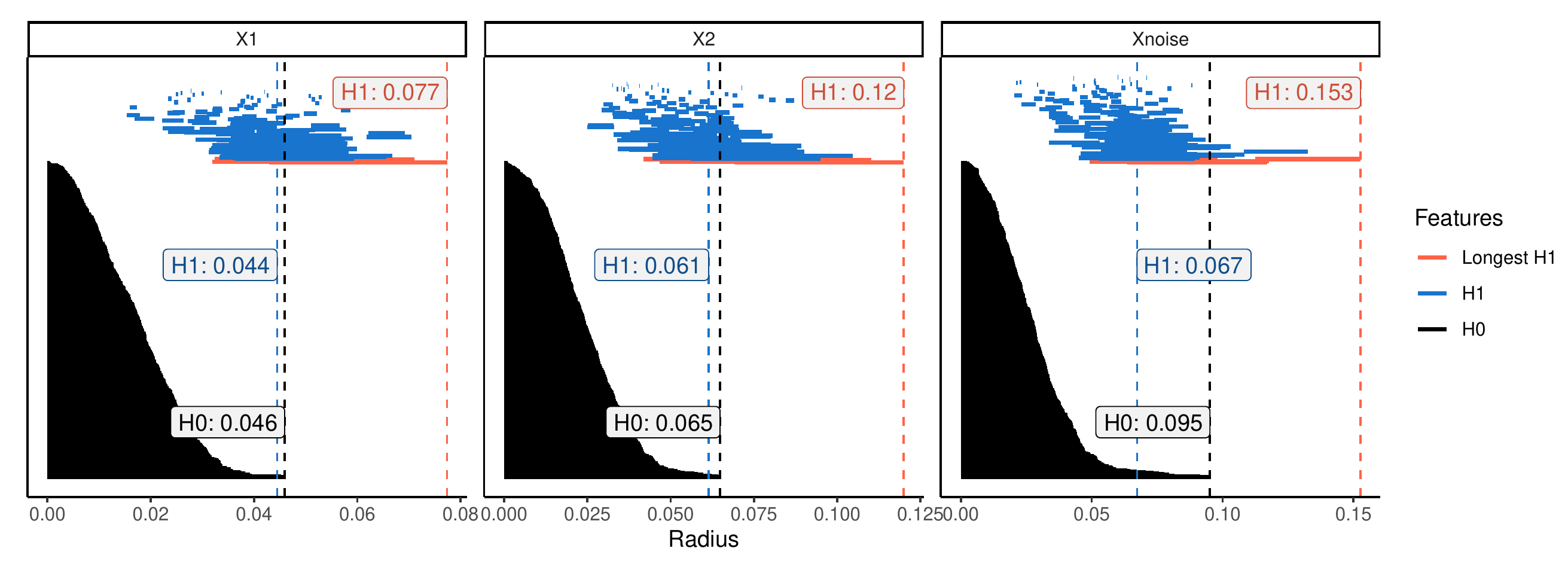}
	\linebreak
	\includegraphics[width=\textwidth]{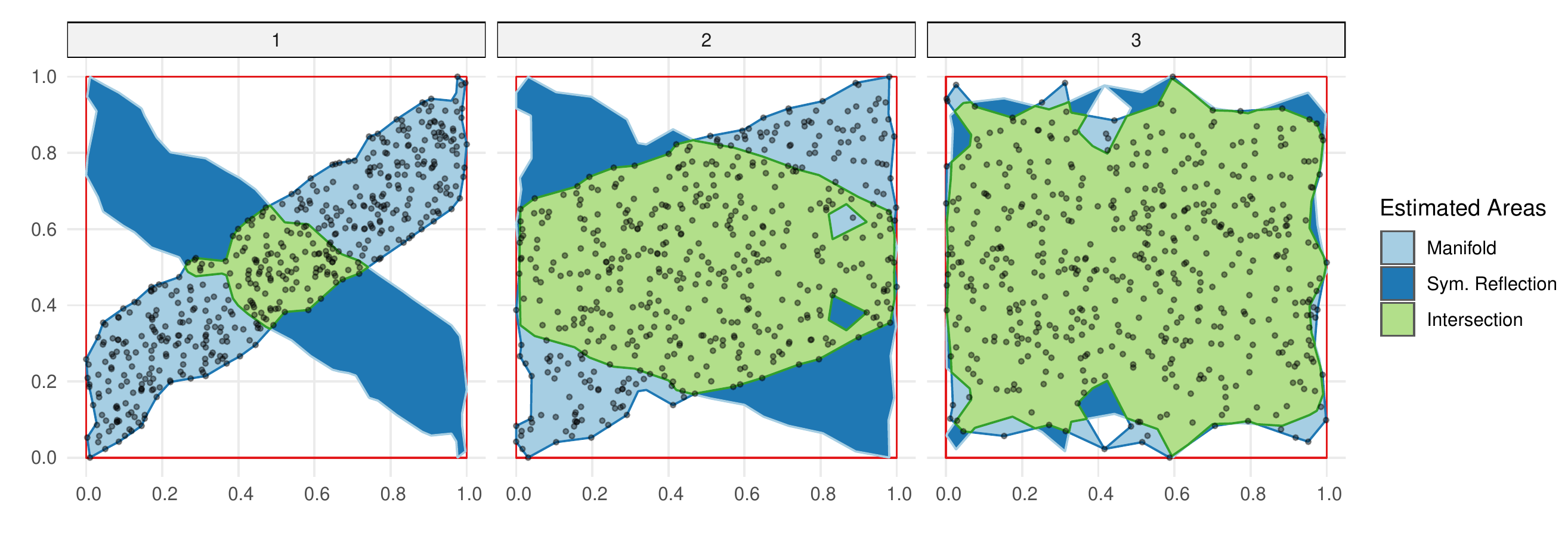}
	\begin{tabular}{cccccc}
		\toprule
		Variable
		& $\varepsilon$
		& $\mathrm{Area}(\mathcal{V}(G))$
		& $\mathrm{Area}(B)$
		& $S_i^{\mathrm{Geom}}$
		& $S_{i}$ \\
		\midrule %
		$X_{1}$ & 0.08 & 0.31 & 1.00 & 0.69 & 0.81 \\
		$X_{2}$ & 0.12 & 0.64 & 1.00 & 0.36 & 0.19 \\
		$X_{3}$ & 0.15 & 0.80 & 0.98 & 0.19 & 0.01 \\
		\bottomrule
	\end{tabular}

	\caption{Results for the Linear case.}
	\label{fig:linear}
\end{figure}

The linear model in Figure~\ref{fig:linear} allow us to identify that
the variable $X_1$ doubles the relevance of $X_{2}$. The indices
$\rho^\mathrm{Geom}$ and $S^\mathrm{Geom}$ are coherent in this case.
The variables $X_{3}$ to $X_{5}$ will have less relevant indices. We
conclude how the empty spaces are present according to the relevance
level of the variable.

\begin{figure}[ht!]
	\centering
	\includegraphics[width=\textwidth]{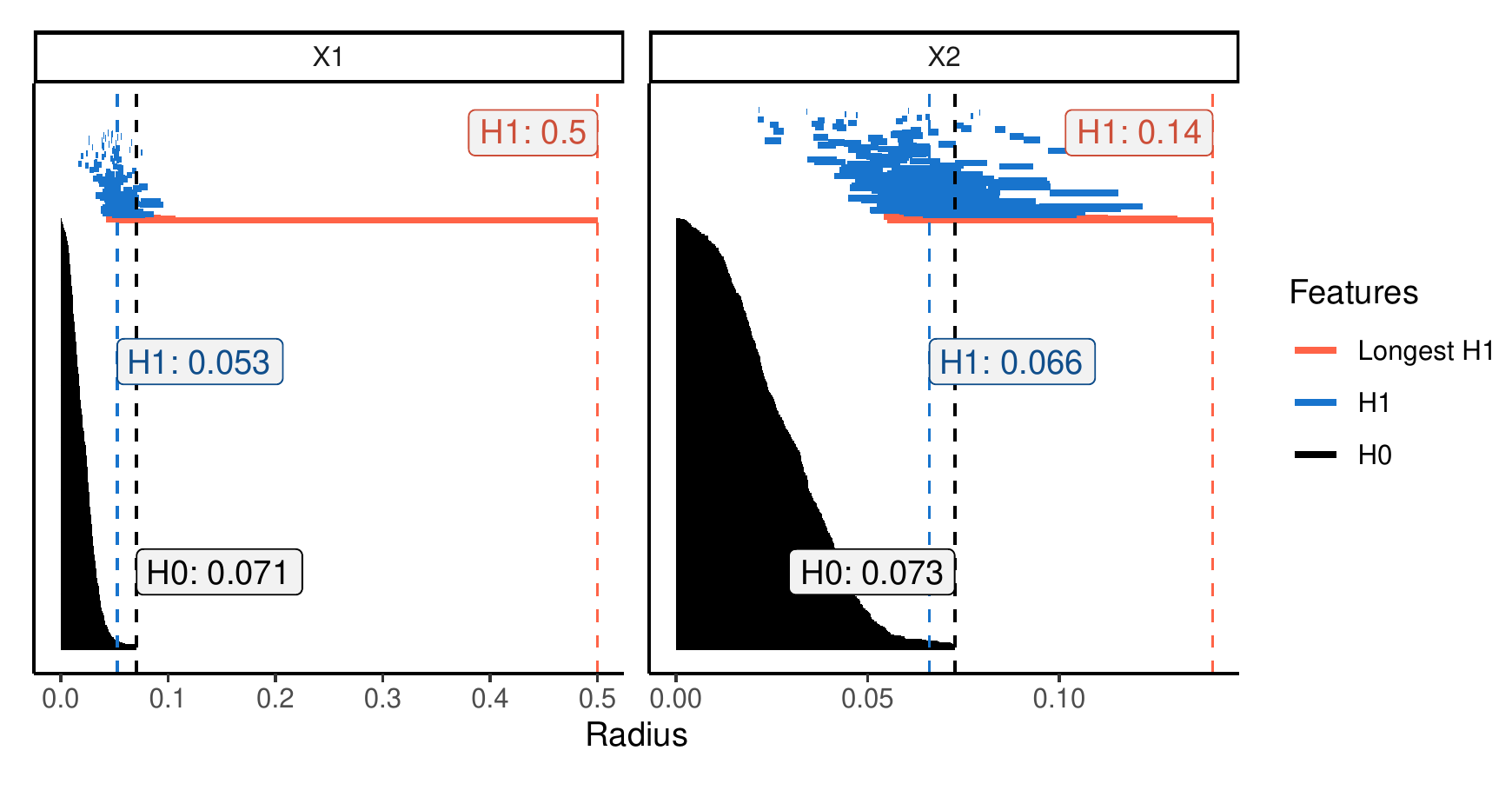}
	\linebreak
	\includegraphics[width=\textwidth]{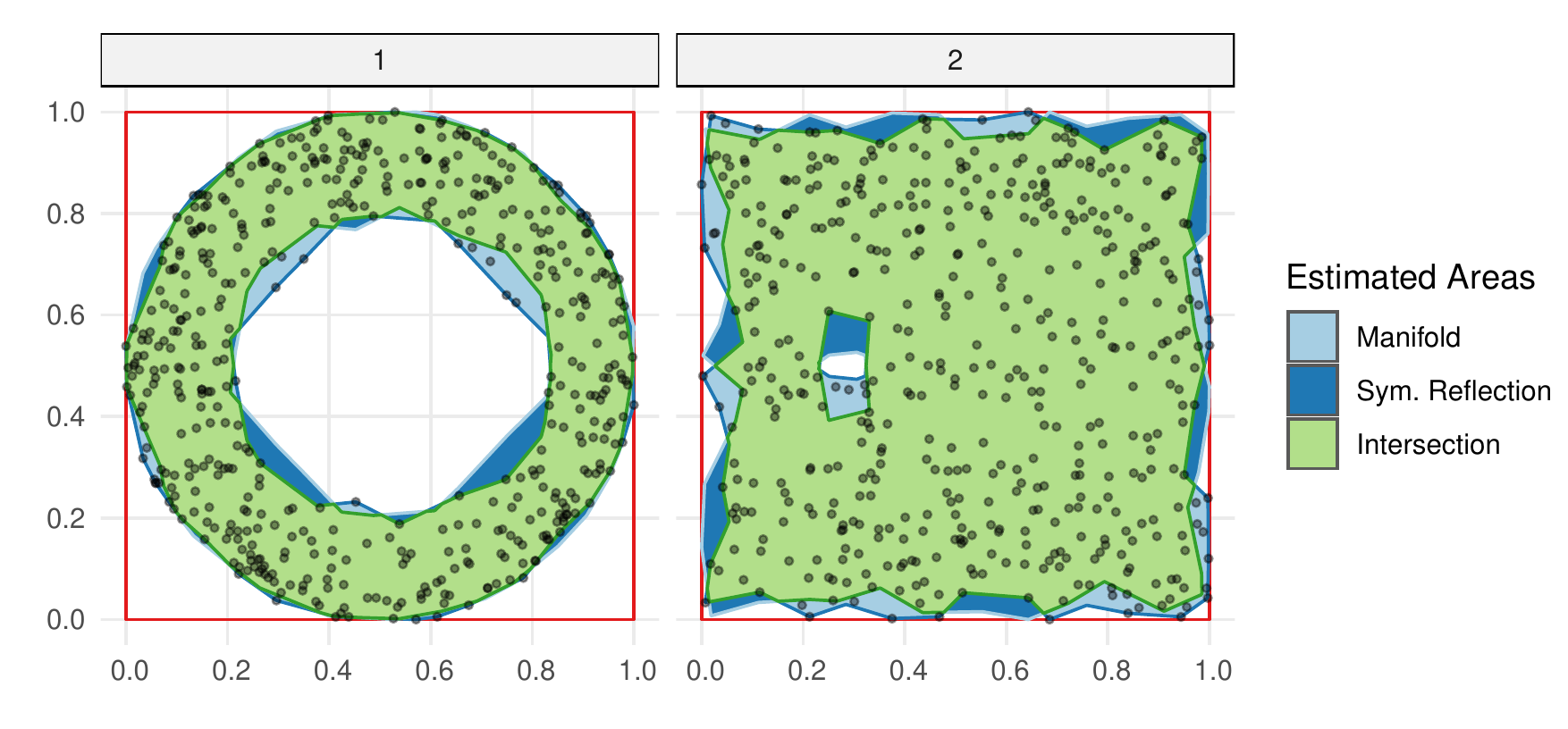}
	\begin{tabular}{cccccc}
		\hline
		Variable
		& $\varepsilon$
		& $\mathrm{Area}(\mathcal{V}(G))$
		& $\mathrm{Area}(B)$
		& $S_i^{\mathrm{Geom}}$
		& $S_{i}$ \\
		\hline
		$X_1$ & 0.10 & 2.07 & 3.94 & 0.48 & 0.03 \\
		$X_2$ & 0.13 & 3.66 & 3.99 & 0.08 & 0.04 \\
		\hline
	\end{tabular}
	\caption{Results for the Circle with hole case.}
	\label{fig:circle}
\end{figure}

The model of a circle with a hole in Figure~\ref{fig:circle} presents
a topological particularity. It is a model were there exist a clear
geometric pattern between $X_1$ and $Y$ but neither $X_1$ nor $X_2$ are
relevant.  Observe how the first variable has $\rho_1^\mathrm{Geom}$
equal to $0.48$ and the $\rho_2^\mathrm{Geom}$ equals to $0.08$.
However, $S_1^\mathrm{Geom}$ and $S_2^\mathrm{Geom}$ are less than
\(0.05\). It means, our estimators detect a geometrical correlation in the
first variable but neither one has impact in the model.

\begin{figure}[ht!]
	\centering
	\includegraphics[width=\textwidth]{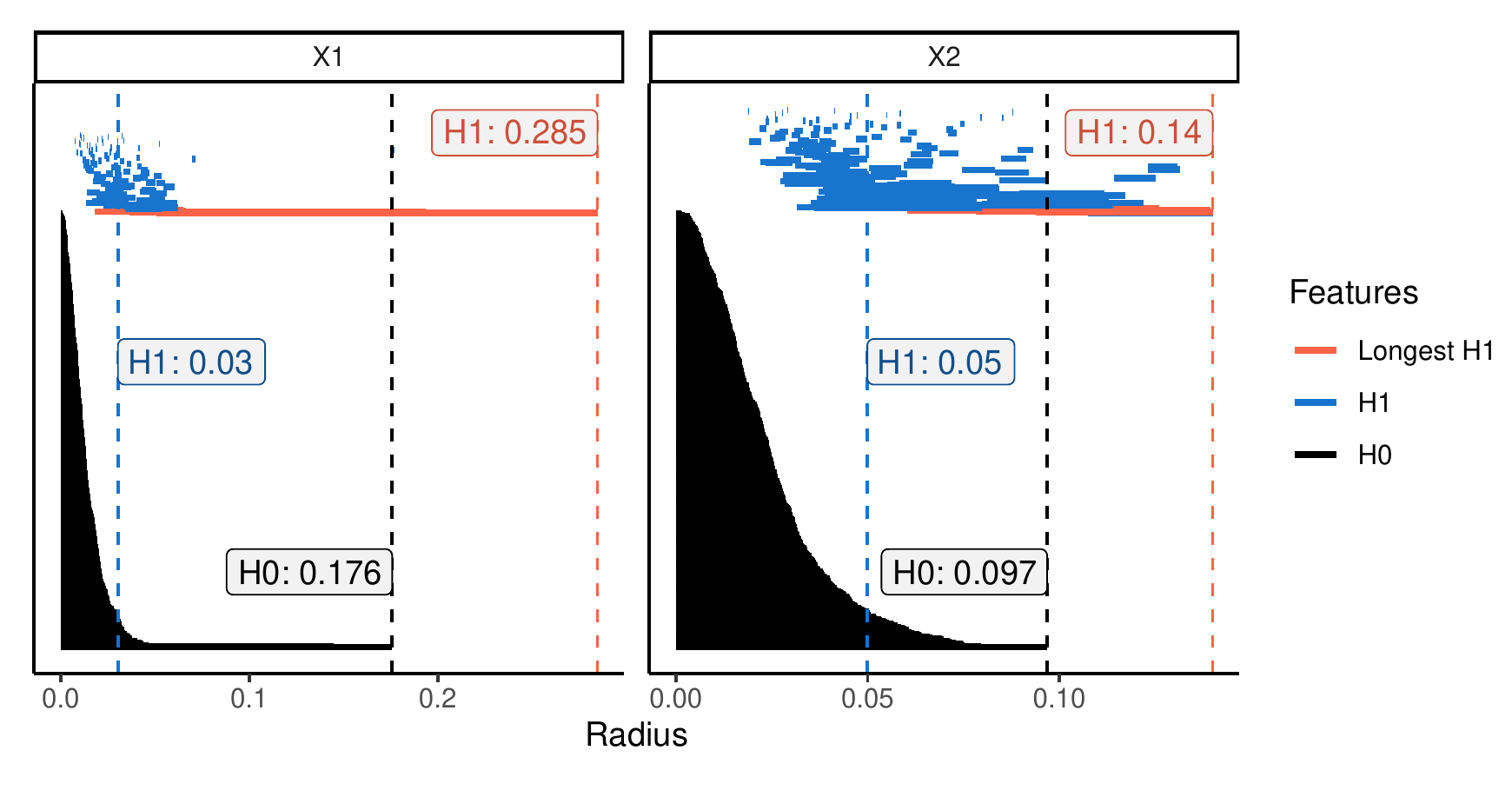}
	\linebreak
	\includegraphics[width=\textwidth]{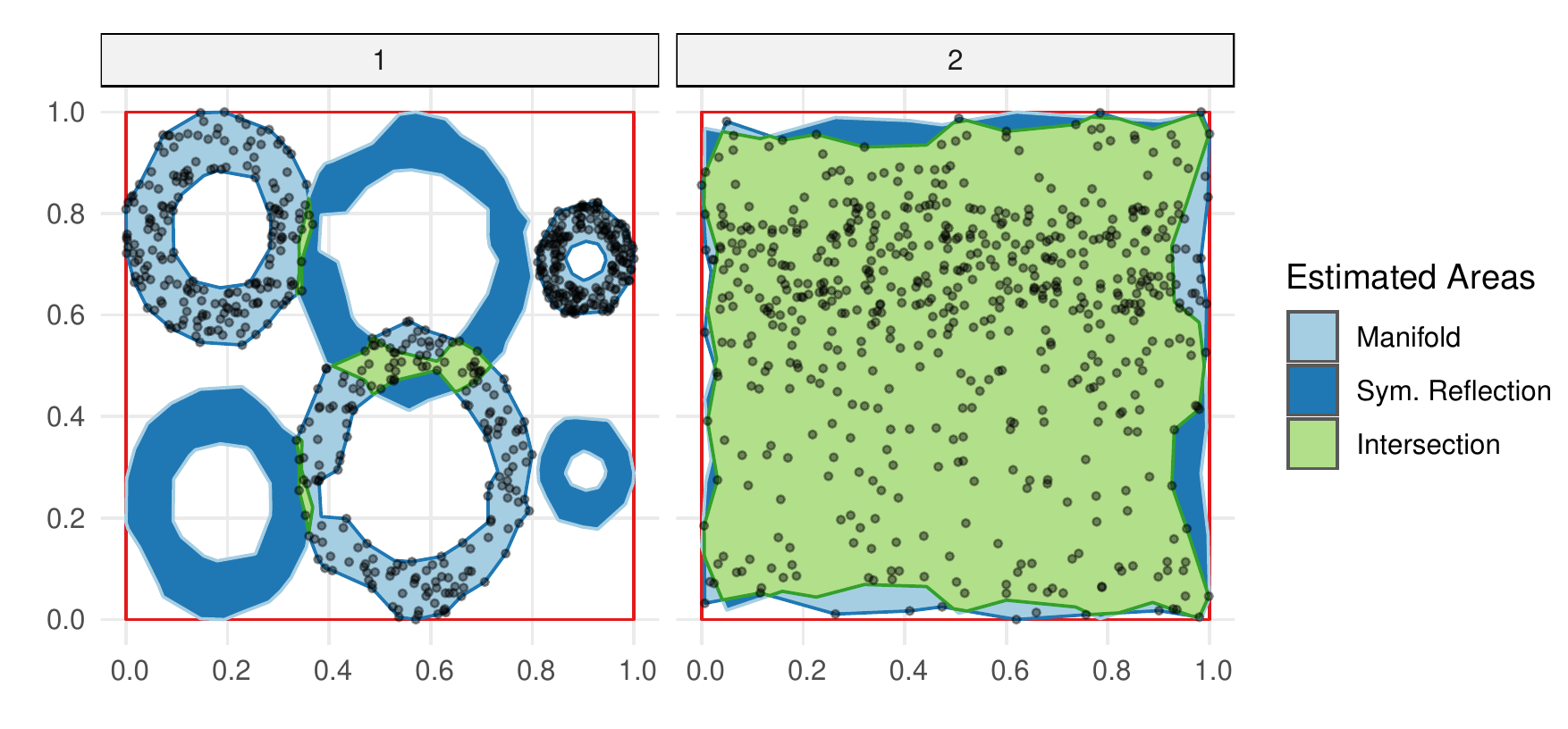}
	\begin{tabular}{cccccc}
		\hline
		Variable
		& $\varepsilon$
		& $\mathrm{Area}(\mathcal{V}(G))$
		& $\mathrm{Area}(B)$
		& $S_i^{\mathrm{Geom}}$
		& $S_{i}$ \\
		\hline
		\(X_{1}\) & 0.09 & 0.23 & 1.00 & 0.77 & 0.64 \\
		\(X_{2}\) & 0.16 & 0.92 & 1.00 & 0.08 & 0.01 \\
		\hline
	\end{tabular}
	\caption{Results for the Circle with two hole case.}
	\label{fig:circle2}
\end{figure}

The other atypical case is a model of connected circles with holes in
Figure~\ref{fig:circle2}. Both circles have different scales and
positions.  In this case, the first variable has a geometric pattern
while the second one is random input. Both the $\rho^\mathrm{Geom}$
and $S^\mathrm{Geom}$ detect that the first variable is the geometric
relevant while the second not.

\begin{figure}[ht!]
	\centering
	\includegraphics[width=\textwidth]{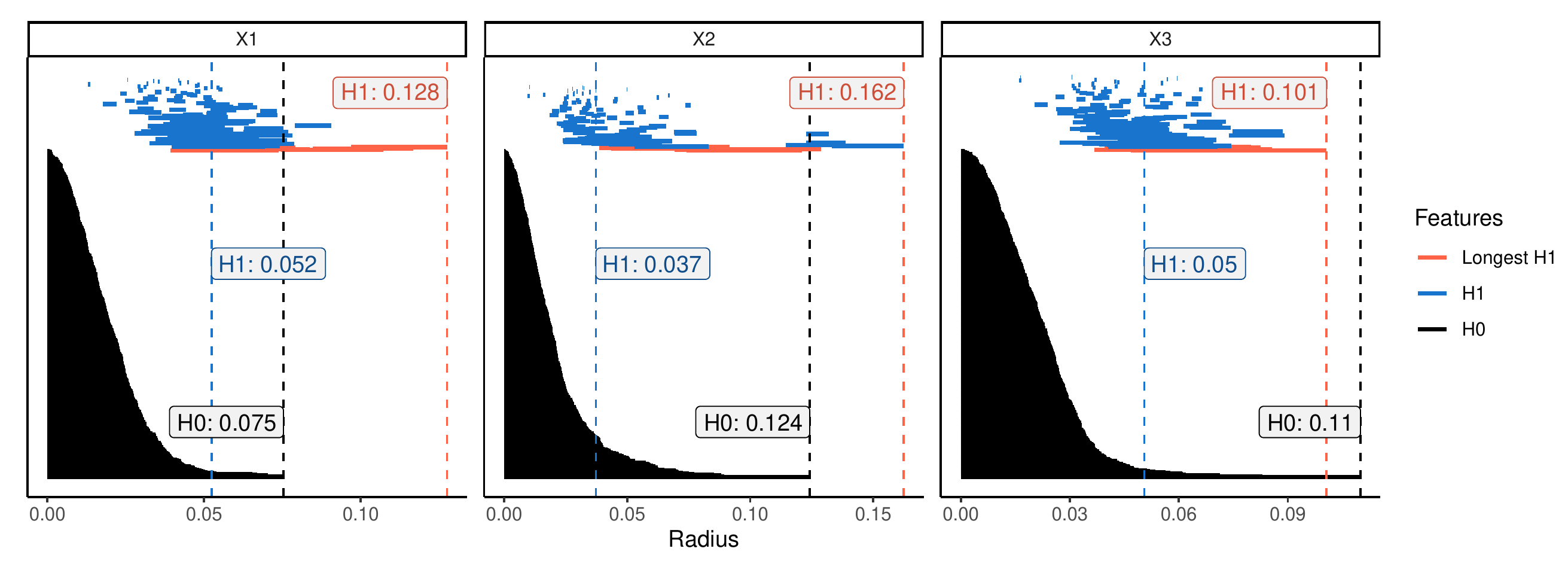}
	\linebreak
	\includegraphics[width=1.2\textwidth]{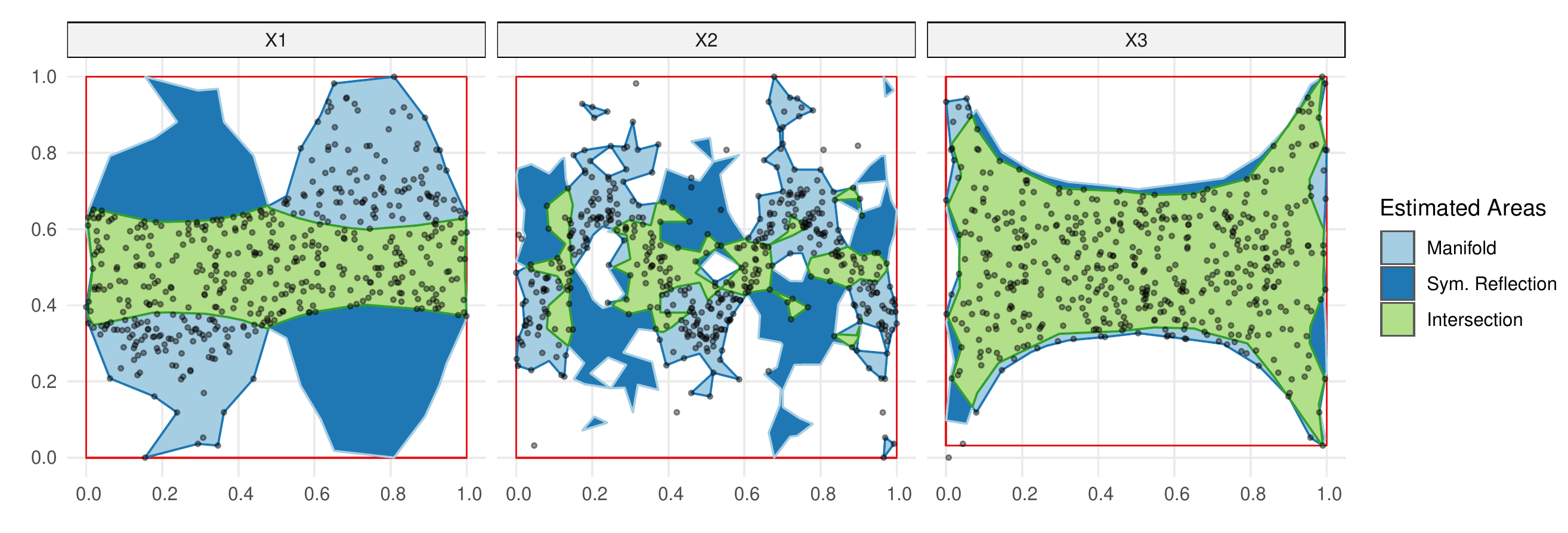}
	\begin{tabular}{cccccc}
		\hline
		Variable
		& $\varepsilon$
		& $\mathrm{Area}(\mathcal{V}(G))$
		& $\mathrm{Area}(B)$
		& $S_i^{\mathrm{Geom}}$
		& $S_{i}$ \\
		\hline
		\(X_{1}\) & 0.16 & 0.48 & 1.00 & 0.49 & 0.29 \\
		\(X_{2}\) & 0.08 & 0.27 & 1.00 & 0.64 & 0.48 \\
		\(X_{3}\) & 0.13 & 0.49 & 0.97 & 0.07 & 0.02 \\
		\hline
	\end{tabular}
	\caption{Results for the Ishigami case.}
	\label{fig:ishigami}
\end{figure}

Another classic model is the Ishigami. model in sensitivity analysis
because it presents a strong non-linearity and non-monotonicity with
interactions in $X_3$. With other sensitivity estimators the variables
$X_{1}$ and $X_{2}$ have great relevance to the model, while the third
one $X_{3}$ has almost zero. For a further explanation of this
function we refer the reader to~\cite{SobolUse1999}.
Figure~\ref{fig:ishigami}. In our case all the geometric correlations
are between  0.5 to 0.7.  As noticed by \cite{HernandezGeometrical2019}, this
measure detects the geometric structure present in every variable.
But, this measure gives an incomplete perspective of the problem
because is not saying if they are relevant or not to the model. For
this, the index $S^\mathrm{Geom}$ establish the order of relevance of
the variables as: $X_2, X_1, X_3$. This match with the theoretical
case discussed at the beginning of Section~\ref{sec:methodology}.

 \FloatBarrier

\section{Conclusions}
\label{sec:conclusions}

In this work we proposed a simple method to detect the relevance of a variable
in a model. We estimated first the embedding manifold in \(\mathbb{R}^{2}\) of
each input variable and the output variable \(Y\). We use the Vietoris-Rips as a
proxy for the intrinsic geometry of the object. By application of simple affine
transformations we compared the area of the estimated manifold to that of its
symmetric reflection and use this information to create a score. If their
symmetric difference were small in terms of area relative to the whole geometric
pattern, the variable's relevance were small as well. Otherwise it would be
higher.

As our experiments show, the geometric sensitivity index constructed in our
method detected the geometric structure within the point cloud. Even if the
construction depends on the simplicial complexes, we based the method on the
area of the complex. This simplified the understanding of the index for
non-technical users. Our index calculated the empty space created by the data
after re-scaling to the the interval \([0,1]\times[0,1]\). This approach allowed
us to normalize the total space taken by the data. If the space is small, then
there is not an obvious pattern between the variable and their response.
Otherwise, we can recognize patterns that influenced the value of \(Y\).

One caveat of this method was the choice of radius in order to construct the
Vietoris-Rips complex. This parameter is the most sensible part of the method,
and yet the most difficult to select. In our examples, we first observed the
barcode and selected the radius according to the most prominent features on the
data. An automated approach to this selection can be done and can be implemented
as part of the algorithm.

While the method provided here showed to be of use, we recognize that its
extension into higher dimension modeling and for the analysis of multiple
variables at a time can be very costly in a computational sense. Another issue
would be to determine the hyperspace symmetrization and the subsequent
calculation of hyper volumes. One plausible approach to solving this problem
would either be projecting the data in low dimensional spaces by means of
principal component analysis, projection pursuit or following an specific
direction in a grand tour (\cite{CookGrand2008}).

Other natural approach to extend this method to a broader class of problems
will follow by exploiting the deeper topological features of the data. Tools like
Euler curves, landscapes or cohomology analysis would presumably lead to new
insights on the relevance of each variable in the model.

\section*{Declarations}

\subsection*{Funding}

First and second authors acknowledge the financial support from
Escuela de Matem\'atica de la Universidad de Costa Rica, through
CIMPA, Centro de Investigaciones en Matem\'atica Pura y Aplicada
through the projects 821--B7--254 and 821--B8--221 respectively.

Third author acknowledges the financial support from Escuela de
Matem\'atica de la Universidad de Costa Rica, through CIMM, Centro de
Investigaciones Matem\'aticas y Metamatem\'aticas through the
project 820--B8--224.

\subsection*{Code availability}

All the calculations in this package were made using an own R-package
\texttt{topsa}. It estimates sensitivity indices reconstructing the embedding
manifold of the data. The reconstruction is done via a Vietoris Rips with a
fixed radius. Then the homology of order 2 and the indices are estimated. The
package is published on CRAN on this address
\url{https://cran.r-project.org/package=topsa}

\subsection*{Conflict of interests}

On behalf of all authors, the corresponding author states that there is no
conflict of interest.

\bibliographystyle{plain} %
\bibliography{references-symmetric-reflection.bib}






\end{document}